\documentclass[a4paper,leqno]{article}

\usepackage{amsmath}
\usepackage{amsthm}
\usepackage{amssymb}
\usepackage{amsfonts}
\usepackage[applemac]{inputenc}
\usepackage{mathrsfs}
\usepackage{color}
\usepackage[hang,flushmargin]{footmisc} 

\def\R{{\mathbb R}}
\def\Z{{\mathbb Z}}
\def\E{{\mathbb E}}

\def\le{\leqslant}
\def\ge{\geqslant}

\DeclareMathOperator*{\osc}{osc}

\DeclareMathOperator{\arsinh}{arsinh}

\theoremstyle{plain}
\newtheorem{theorem}{Theorem}[section]
\newtheorem{lemma}[theorem]{Lemma}

\newtheorem{proposition}[theorem]{Proposition}

\theoremstyle{definition}

\newtheorem*{remark*}{Remark}

\numberwithin{equation}{section}

\begin{document}

\title{On annealed elliptic Green function estimates}
\author{Daniel Marahrens\thanks{Max Planck Institute for Mathematics in the Sciences, Inselstr.\ 22, 04103 Leipzig, Germany, \texttt{Daniel.Marahrens@mis.mpg.de} resp.\ \texttt{Felix.Otto@mis.mpg.de}} \and Felix Otto\footnotemark[1]}

\date{\today}

\renewcommand{\thefootnote}{\fnsymbol{footnote}} 
\footnotetext{\noindent\today\\\emph{MSC2010 subject classifications.} Primary: 35B27, Secondary: 35J08, 39A70, 60H25}     
\renewcommand{\thefootnote}{\arabic{footnote}}

\maketitle





\begin{abstract} 
 We consider a random, uniformly elliptic coefficient field $a$ on the lattice $\Z^d$. The distribution $\langle \cdot \rangle$ of the coefficient
field is assumed to be stationary. Delmotte and Deuschel showed that the gradient and second mixed derivative of the parabolic Green
function $G(t,x,y)$ satisfy optimal annealed estimates which are $L^2$ resp.\ $L^1$ in probability, i.e.\ they
obtained bounds on $\langle |\nabla_x G(t,x,y)|^2 \rangle^{\frac{1}{2}}$ and $\langle |\nabla_x \nabla_y G(t,x,y)| \rangle$, see  T.\ Delmotte and J.-D.\ Deuschel: On estimating the derivatives of
symmetric diffusions in stationary random environments, with applications to the $\nabla\phi$ interface model, Probab.\ Theory Relat.\ Fields {\bf 133} (2005), 358--390. In particular,
the elliptic Green function $G(x,y)$ satisfies optimal annealed bounds. In a recent work, the authors extended these elliptic bounds to higher
moments, i.e.\ $L^p$ in probability for all $p<\infty$, see D.\ Marahrens and F.\ Otto: {Annealed estimates on the Green function}, arXiv:1304.4408 (2013). In this note, we present a new argument that relies purely on elliptic theory to
derive the elliptic estimates (see Proposition \ref{prop:DD_ell} below) for $\langle |\nabla_x G(x,y)|^2 \rangle^{\frac{1}{2}}$ and $\langle
|\nabla_x \nabla_y G(x,y)| \rangle$.
\end{abstract}

\section{Introduction}
In this work, we consider linear second-order difference equations with uniformly elliptic, random coefficients of the form
\begin{equation}\label{ell}
 \nabla^*(a \nabla u)(x) = f(x) \quad \text{for all $x\in\Z^d$.}
\end{equation}
If there is no danger of confusion, we also write $\nabla^* a \nabla u$ for $\nabla^*(a\nabla u)$. 
In this discrete difference equation, the gradient $\nabla$ and the (negative) divergence $\nabla^*$ on $\Z^d$ are defined as follows: Let
$\E^d$ denote the set of \emph{edges} of $\Z^d$ consisting of all pairs $[x,x+e_i]$ of neighboring vertices with $x\in\Z^d$,
$i=1,\ldots,d$. Here $e_1, \ldots, e_d$ is the canonical basis of $\R^d$. Then we set
\begin{align*}
 \nabla \zeta ([x,x+e_i]) &= \zeta(x+e_i)-\zeta(x),\\
 \nabla^*\xi(x) &= \sum_{i=1}^d \big(\xi([x-e_i,x]) - \xi([x,x+e_i])\big),
\end{align*}
for all fields on vertices $\zeta : \Z^d \to \R$ (which we think of as scalar fields) and fields on edges $\xi : \E^d \to \R$ (which we think of as vector fields).
In general, we will denote edges by the letters $e$ and $b$ and vertices by the letters $x$, $y$ and $z$. The operators
$\nabla$ and $\nabla^*$ are adjoint in the sense of
\[
 \sum_{e\in\E^d} \xi(e) \nabla\zeta(e) = \sum_{x\in\Z^d} \zeta(x) \nabla^*\xi(x).
\]
In \eqref{ell}, the coefficient field $a$ is a field on edges $a:\E^d \to \R$.

\smallskip

Our assumption on the coefficient field is two-fold: one deterministic and one probabilistic assumption. The deterministic assumption is
that of \emph{uniform ellipticity}:
We assume that for every $e\in\E^d$ we have that $\lambda \le a(e) \le 1$. Here $\lambda\in(0,1)$ denotes the ellipticity ratio
which is fixed  throughout the paper. We denote the space of uniformly elliptic coefficient fields $a:\E^d \to [\lambda,1]$ by
$\Omega$, i.e.\ we set
\[\label{Omega}
 \Omega:= [\lambda,1]^{\E^d}.
\]
In this work, the coefficient field is assumed to be distributed according to a probability measure on $\Omega$. Following the convention in
statistical mechanics, we call this probability measure an \emph{ensemble} and denote its expectation by $\langle \cdot \rangle$.
The probabilistic assumption on $\langle \cdot \rangle$ is that of \emph{stationarity}. To define this property, we note that $\Z^d$ acts on
$\E^d$ by
translation and we denote by $e+x\in\E^d$ the edge $e\in\E^d$ shifted by $x\in\Z^d$. Then stationarity means
that the coefficient field is distributed according to some probability measure on $\Omega$ such that $a$ and
$a(\cdot+x)$ have the same distribution for all $x\in\Z^d$.

\smallskip


We are interested in proving Green function estimates. The Green function $G(a;x,y)=G(x,y)$ is the fundamental solution of
\eqref{ell}, i.e.\ the solution to
\begin{equation}\label{Green}
 \nabla^*(a \nabla G(\cdot,y))(x) = \delta(x-y),
\end{equation}
where the right hand side is the discrete Dirac on $\Z^d$ defined as
\[
 \delta(x) = \begin{cases}
              1 & x=0,\\
	      0 &\text{otherwise.}
             \end{cases}
\]
Dimension $d=2$ needs a bit more care in terms of the definition of the Green function. Since in this work we are only interested in
\emph{gradient} estimates, this is merely technical and will be ignored here. 
It is well-known since the work of Nash \cite{Nash} and Aronson \cite{Aro} that in dimension $d>2$ the Green function $G$ itself
satisfies
\begin{equation}\label{nash}
 \frac{1}{C} (|x-y|+1)^{2-d} \le G(x,y) \le C (|x-y|+1)^{2-d},
\end{equation}
for some constant $C=C(d,\lambda)$ depending only on the dimension $d$ and the ellipticity contrast $\lambda$. Here and throughout, we
denote (generic) constants that only depend on their arguments $(\cdot)$ by $C(\cdot)$. In particular, the bounds \eqref{nash} are
\emph{quenched} bounds, i.e.\ they do not depend on the choice of $a\in\Omega$. The bounds are optimal since they are the same as the bounds
for the
constant coefficient Green function. On the other hand, without further assumptions on the coefficients besides uniform ellipticity,
we cannot expect the same bounds as for the constant coefficient Green function to hold for the \emph{gradients} of non-constant
coefficient Green functions uniformly in $a\in\Omega$. In fact, de Giorgi-Nash-Moser theory only
yields 
\begin{align}\label{degiorgi}
 |\nabla G(e,y)| &\le C(d,\lambda) (|e-y|+1)^{2-d-\alpha_0},\\
 |\nabla \nabla G(e,b)| &\le C(d,\lambda) (|e-b|+1)^{2-d-2\alpha_0},\label{degiorgi2}
\end{align}
for some $\alpha_0=\alpha_0(d,\lambda)>0$. Let us explain the notation here: Since $G$ is a function of two variables, we make the
convention that in case of ambiguity we always let the derivative fall onto the edge
variable. For instance, in $\nabla G(e,y)$ the derivative is taken in the first variable along the edge $e\in\E^d$. The term $\nabla\nabla
G(e,b)$ denotes the second \emph{mixed} derivative, i.e.\ we take one
derivative in the first variable of $G$ along the edge $e$ and one derivative in the second variable of $G$ along the edge $b$.
By an abuse of notation, $|e-y|$ denotes the distance of the origin from the mid-point of $e-y\in\E^d$ and $|e-b|$ denotes the
distance between the two mid-points of the edges $e,b\in\E^d$. We cannot expect more than \eqref{degiorgi} and \eqref{degiorgi2},
since the constant coefficient 
bounds would imply in particular almost Lipschitz-continuity of $a$-harmonic functions, cf.\ Corollary 4 in \cite{MO}. This is where
stationarity comes into play. Indeed, for the parabolic Green function, i.e.\ the solution to
\begin{equation}\label{Green_para}
 \partial_t G(t,x,y) + (\nabla^* a \nabla G(t,\cdot,y))(x) = 0, \qquad \text{s.t.} \quad G(0,x,y) = \delta(x-y),
\end{equation}
Delmotte and Deuschel \cite{DD} have shown \emph{annealed} (i.e.\ in mean) Green function
estimates, which are the content of the following Proposition.
\begin{proposition}[Theorem 1.1 of Delmotte and Deuschel \cite{DD}]\
 If the ensemble $\langle \cdot \rangle$ on $\Omega$ is stationary, then we have that
\begin{align}\label{DD}
 \langle |\nabla G(t,e,y)|^{2} \rangle^{\frac{1}{2}} &\le C(d,\lambda) (t+1)^{-\frac{d-1}{2}} \exp\big( {-\textstyle \frac{1}{C(d,\lambda)}
\eta(t,|e-y|)} \big) \quad\text{and}\\
 \langle |\nabla \nabla G(t,e,b)| \rangle &\le C(d,\lambda) (t+1)^{-\frac{d-2}{2}} \exp\big( {-\textstyle \frac{1}{C(d,\lambda)}
\eta(t,|e-b|)}
\big).\label{DD2}
\end{align} 
for all $t>0$, all $e,b\in\E^d$ and all $y\in\Z^d$, where $\eta(t,r) = r \arsinh (\frac{r}{t}) - ( \sqrt{t^2 + r^2} - t^2 )$.
\end{proposition}
Thus the behaviour of the parabolic Green function is slightly more complicated than in the continuum setting, where it is simply Gaussian. This phenomenon is due to discrete (finite size) effects, which for $|x| \gg t$ allow the Green function to spread out much faster than expected. Indeed, the function $\exp (-\eta(t,r))$ behaves like $\exp(-r^2/t)$ for small $r/t$ and like $(t/r)^r \exp r$ for large $r/t$, cf.~\cite[Remark 2]{Del}.
Note that these estimates require only uniform ellipticity and stationarity of the coefficient field $a$. Since $G(x,y) = \int_0^\infty
G(t,x,y) \;dt$, the estimates \eqref{DD} and \eqref{DD2} immediately imply annealed bounds on the elliptic Green function, for which we
aim in this note to provide an alternative, self-contained and direct proof.
\begin{proposition}\label{prop:DD_ell}
 If the ensemble $\langle \cdot \rangle$ on $\Omega$ is stationary, then we have that
\begin{align}\label{DD_ell}
 \langle |\nabla G(e,y)|^{2} \rangle^{\frac{1}{2}} &\le C(d,\lambda) (|e-y|+1)^{1-d} \quad\text{and}\\
 \langle |\nabla \nabla G(e,b)| \rangle &\le C(d,\lambda) (|e-b|+1)^{-d}.\label{DD_ell2}
\end{align}
for all $e,b\in\E^d$ and $y\in\Z^d$.
\end{proposition}
In \cite{MO}, the authors upgraded \eqref{DD_ell} and \eqref{DD_ell2} to higher moments under under the assumption of a logarithmic Sobolev
inequality (LSI) with constant $\rho$. The ensemble $\langle \cdot \rangle$ is said to satisfy a LSI with constant $\rho > 0$ if
\begin{equation}\label{LSI}
 \Big\langle \zeta^2 \log {\textstyle\frac{\zeta^2}{\langle \zeta^2 \rangle}} \Big\rangle \le \frac{2}{\rho} \Big\langle \sum_{e\in\E^d}
\Big( \osc_{a(e)} \zeta \Big)^2 \Big\rangle
\end{equation}
for all random variables $\zeta : \Omega \to \R$. Here the oscillation is taken over all coefficient fields
$\{\tilde a(b)\}_{b\in\E^d}\in\Omega$ that coincide with $a$ outside of $e$, i.e.\
\begin{multline*}
 \osc_{a(e)} \zeta(a) := \sup \{ \zeta(\tilde a) : \tilde a \in \Omega \text{ s.\ t.\ } \tilde a(b) = a(b) \text{ for all } b\neq e \}\\ -
\inf \{ \zeta(\tilde a) : \tilde a \in \Omega \text{ s.\ t.\ } \tilde a(b) = a(b) \text{ for all } b\neq e \}.
\end{multline*}
As was shown in \cite{MO}, the LSI is satisfied for all
identically and independently distributed random coefficient fields. The main result of \cite{MO} is the following.
\begin{proposition}[Theorem 1 in \cite{MO}]\label{prop:MO}
Let the ensemble $\langle \cdot \rangle$ on $\Omega$ be stationary and satisfy \eqref{LSI}. Then we have that
\begin{align*}
 \langle |\nabla G(e,y)|^{2p} \rangle^{\frac{1}{2p}} &\le C(d,\lambda,\rho,p) (|e-y|+1)^{1-d} \quad\text{and}\\
 \langle |\nabla \nabla G(e,b)|^{2p} \rangle^{\frac{1}{2p}} &\le C(d,\lambda,\rho,p) (|e-b|+1)^{-d}
\end{align*}
for all $x,y\in\Z^d$ and all $p<\infty$.
\end{proposition}
As remarked before, the elliptic estimates \eqref{DD_ell} and \eqref{DD_ell2} follow from integrating the parabolic estimates \eqref{DD}
and \eqref{DD2} of Delmotte and Deuschel \cite{DD}. Indeed, the parabolic estimates are the only point in \cite{MO} where parabolic theory
enters. In this note, we point out an alternative approach to obtain Proposition~\ref{prop:DD_ell}, i.e.\ \eqref{DD_ell} and
\eqref{DD_ell2}, based purely on elliptic theory. This has clear conceptual advantages but we also believe that the proof in itself may be
interesting to the reader.
The proof of Proposition \ref{prop:DD_ell} relies on the following quenched result, showing optimal spatially averaged decay on annuli.
\begin{lemma}\label{lem:quenched}
For all $a\in\Omega$, vertices $y\in\mathbb{Z}^d$ and radii $R\ge 1$, we have that
\begin{align}
\bigg(R^{-d}\sum_{e:R\le |e-y|\le 2R}|\nabla G(e,y)|^2\bigg)^\frac{1}{2}&\le C(d,\lambda) R^{1-d},\label{quenched}\\
\bigg(R^{-2d}\sum_{e:8R\le|e-y|\le 16R}\sum_{b:|b-y|\le R}|\nabla \nabla G(e,b)|^2\bigg)^\frac{1}{2}&\le C(d,\lambda) R^{-d}.
\label{quenched2}
\end{align}
These estimates are optimal by comparison with the constant-coefficient case.
%
\end{lemma}
This lemma is inspired by the work of the second author, Lamacz and Neukamm \cite{LNO} on degenerate equations related to percolation and we
shall prove it in Section \ref{sec:quenched}. In Section \ref{sec:DD_ell}, we will use stationarity to deduce
Proposition~\ref{prop:DD_ell} from Lemma~\ref{lem:quenched}.

%

\section{Proof of Lemma \ref{lem:quenched}}\label{sec:quenched}

For simplicity, we will first establish the result in the continuum case by standard
arguments in the spirit of De Giorgi. The sole additional difficulty coming from discreteness is the absence of the Leibniz rule and the chain rule. We will not worry about regularity and finiteness
in the continuum setting and address this only when indicating the necessary changes for the discrete case.
Only the proof of Step 2 produces additional lower order terms in the discrete case because of the cut-off function, cf.\ Step 3. In both the discrete and continuum treatment, we follow \cite{LNO}.
The symbols $\nabla$ and $\nabla^*$ now denote the continuum gradient and its formal transpose, the continuum (negative) divergence.
Furthermore, $\lesssim$ stands for a generic constant that only depends on $d$ and $\lambda$.
We fix an arbitrary coefficient field $a\in\Omega$.


\medskip

{\bf Step 1}. Suppose that we are given $u$ and $f$ such that
\begin{equation}\label{S.4}
\nabla^*a\nabla u=f\quad\mbox{in}\;\mathbb{R}^d.
\end{equation}
Then we have for any ball $B_R$ of radius $R$:
\begin{equation}\label{S.20}
 \inf_{\bar u \in \R}\int_{B_R}|u - \bar u| \lesssim R^2\int_{\mathbb{R}^d} |f| .
\end{equation}
Without loss of generality, we may assume by homogeneity that
\begin{equation}\label{S.5}
\int_{\mathbb{R}^d}|f|=1 .
\end{equation}

By modifying $u$ by an additive constant, we may assume that the median of $u$ on $B_R$ vanishes, that is,
\begin{equation}\label{S.6}
|\{u> 0\}\cap B_R|,|\{u<0\}\cap B_R|\le\frac{1}{2}|B_R|.
\end{equation}
Hence it is now enough to show that
\begin{equation}\label{S.10}
\int_{B_R} |u| \lesssim R^2.
\end{equation}
%

\smallskip

After these preparations, we may now start with the actual argument. For given $M\in(0,\infty)$, we test
(\ref{S.4}) with $v_M=\min\{\max\{u,0\},M\}$. Then $v_M$ satisfies
\[
 \int_{\R^d} \nabla v_M\cdot a\nabla v_M  = \int_{\R^d} \nabla u \cdot a\nabla v_M .
\]
Integration by parts yields
\begin{equation}\label{int_by_parts}
\int_{\R^d}\nabla v_M\cdot a\nabla v_M =\int_{\mathbb{R}^d}v_M \nabla^* a\nabla u  = \int_{\mathbb{R}^d} v_M f.
\end{equation}
(Incidentally, the construction of $v_M$ is the only place in this paper where we rely on the fact that we deal
with a \emph{scalar} equation as opposed to an elliptic system.) Note that it is not a priori clear that this integration by parts is
valid; we will justify it in Step 4. 
Using the uniform ellipticity on the l.\ h.\ s.\ and (\ref{S.5}) on the r.\ h.\ s., we obtain the estimate
\begin{equation}\label{S.7}
\lambda\int_{B_R}|\nabla v_M|^2  \le \lambda\int_{\mathbb{R}^d}|\nabla v_M|^2 \le M.
\end{equation}
%
Letting $\overline{v}_M$ denote the average of $v_M$ over $B_R$, we obtain by a Poincar\'e-Sobolev estimate
\begin{equation}\label{S.8}
\bigg( R^{-d} \int_{B_R}|v_M - \overline{v}_M|^{2q} \bigg)^\frac{1}{2q} \lesssim \bigg( R^{-d} \int_{B_R}|R \nabla v_M|^2 \bigg)^\frac{1}{2}
\end{equation}
for some exponent $q>1$ that only depends on $d$. Such a Poincar\'e-Sobolev estimate holds for all $1\le q \le \frac{d}{d-2}$ if $d>2$ and
all $q<\infty$ if $d=2$.
Because of (\ref{S.6}), we have $|\{v_M=0\}\cap B_R|\ge\frac{1}{2}|B_R|$ and hence $\overline{v}_M \le
\frac{M}{2}$. It follows that
\begin{equation}\label{mean_v}
M-\overline{v}_M \ge \frac{M}{2}.
\end{equation}
Finally, we have by definition of $v_M$ and Chebyshev's inequality
\[
 | \{u\ge M \}\cap B_R | \le (M - \overline{v}_M)^{-2q} \int_{B_R} | v_M - \overline{v}_M |^{2q} .
\]
Now \eqref{mean_v} and division by $R^{-d}$ yields
\begin{equation}\label{S.9}
 R^{-d} | \{u\ge M \}\cap B_R | \lesssim M^{-2q} R^{-d} \int_{B_R} | v_M - \overline{v}_M |^{2q} .
\end{equation}
Into \eqref{S.9}, we insert first \eqref{S.8} and then \eqref{S.7} to obtain
\begin{equation}\nonumber
 R^{-d} |\{u\ge M\}\cap B_R| \lesssim M^{-q} R^{(2-d)q},
\end{equation}
which because of $|B_R|\lesssim R^d$ upgrades to
\begin{equation}\nonumber
|\{u\ge M\}\cap B_R|\lesssim R^d \min\{M^{-q} R^{(2-d)q}, 1\}.
\end{equation}
Integration of this estimate in $M\in(0,\infty)$ yields a bound for $u_+:=\max\{u,0\}$, the positive part of $u$:
\begin{equation}\nonumber
\int_{B_R}u_+=\int_0^\infty|\{u\ge M'\}\cap B_R|dM'\lesssim R^d \bigg( \int_M^\infty (M')^{-q} dM' R^{(2-d)q}  + M \bigg)
\end{equation}
for every $M>0$.
Since $q>1$, we may compute the integral and set $M=R^{2-d}$ to obtain
\[
 \int_{B_R}u_+ \lesssim R^d \Big( M^{1-q} R^{(2-d)q} + M \Big) \lesssim R^2.
\]
Symmetry w.\ r.\ t.\ interchange of $u$ and $-u$ yields \eqref{S.20}. 

\smallskip

Changes to the discrete setting: All estimates in this step hold verbatim in the discrete setting. Next to the justification of the
integration by parts, all we need is a discrete version of the Poincar\'{e}-Sobolev estimate \eqref{S.8}. One way to achieve this is to replace balls by boxes and to obtain the discrete version of \eqref{S.8} from the continuum one by applying the latter to piecewise linear, continuous interpolation of the discrete function on a triangulation subordinate to the box in $\Z^d$. Note that even though the discrete and continuum mean values might not coincide, the discrete version of \eqref{S.8} with the continuum mean value (of the interpolation function) implies the one with the discrete mean value by Jensen's inequality. We also remark that since the ball $B_R$ is contained within a box of radius $R$, the result \eqref{S.20} remains valid in the discrete case.

\medskip

{\bf Step 2}. Suppose now that for some ball $B_{2R}$ of radius $2R$:
\begin{equation}\label{S.11}
\nabla^*a\nabla u=0\quad\mbox{in}\;B_{2R}.
\end{equation}
Then we have
\begin{equation}\label{S.17}
\bigg(R^{-d}\int_{B_R}|R\nabla u|^2 \bigg)^\frac{1}{2}\lesssim R^{-d}\int_{B_{2R}}|u| 
\lesssim \bigg(R^{-d}\int_{B_{2R}}|u|^2 \bigg)^\frac{1}{2},
\end{equation}
%
where $B_R$ is the concentric ball of half the radius.
The second estimate follows immediately by Jensen's inequality and hence we just need to prove the first estimate. 
Let $\eta$ denote a cut-off function for $B_R$ in $B_{2R}$, to be further specified below. We test (\ref{S.11})
with $\eta^2u$. Because of the identity
\begin{equation}\label{nk.2}
\nabla(\eta^2u)\cdot a\nabla u=\nabla(\eta u)\cdot a\nabla(\eta u)-u^2\nabla\eta\cdot a\nabla\eta,
\end{equation}
which relies on symmetry of $a$ and that by uniform ellipticity turns into the inequality
\begin{equation}\nonumber
\nabla(\eta^2u)\cdot a\nabla u\ge\lambda|\nabla(\eta u)|^2-u^2|\nabla\eta|^2,
\end{equation}
we obtain
%
%
%
%
%
%
%
\begin{equation}\label{S.12}
\int_{B_{2R}}|\nabla(\eta u)|^2\lesssim \int_{B_{2R}}|\nabla\eta|^2 u^2.
\end{equation}
Estimate (\ref{S.12}) is the standard Caccioppoli estimate.
On the left hand side of (\ref{S.12}), we apply the Sobolev-Poincar\'e inequality 
(this time, on the whole space for functions $v=\eta u$ supported on $B_{2R}$ but with the same
exponent $q>1$ as in Step 1) in form of
\begin{equation}\label{S.14}
\bigg(R^{-d} \int_{\R^d}(\eta u)^{2q}\bigg)^\frac{1}{2q}\lesssim\bigg( R^{-d} \int_{\R^d}|R \nabla(\eta u)|^2\bigg)^\frac{1}{2}.
\end{equation}
On the right hand side, we apply H\"older's inequality (using $q\ge 1$) in form of
\begin{equation}\label{S.13}
\int_{B_{2R}}|R\nabla\eta|^2 u^2
\lesssim 
\bigg(\int_{B_{2R}}\big(|R\nabla\eta|^\frac{2q-1}{q}|u|\big)^{2q}\bigg)^\frac{1}{2q-1}
\bigg(\int_{B_{2R}}                            |u|      \bigg)^\frac{2(q-1)}{2q-1}.
\end{equation}
Since we have the strict inequality $q>1$ and thus $\frac{2q-1}{q}>1$, we can select the cut-off function
such that
\begin{equation}\label{neukamm}
|R\nabla\eta|^\frac{2q-1}{q}\lesssim\eta,
\end{equation}
(for instance, one may select a standard cut-off function $\zeta$ with $|\nabla\zeta|\lesssim 1$ and
then set $\eta=\zeta^r$ with $\frac{1}{r-1}+\frac{1}{q}=1$ so that 
$|R\nabla\eta|\lesssim\eta^\frac{r-1}{r}$ which can be rewritten as \eqref{neukamm})
so that (\ref{S.13}) turns into
\begin{equation}\label{S.15}
R^{-d}\int_{B_{2R}}|R\nabla\eta|^2 u^2
\lesssim 
\bigg(R^{-d}\int_{B_{2R}}|\eta u|^{2q}\bigg)^\frac{1}{2q-1}
\bigg(R^{-d}\int_{B_{2R}}|u|          \bigg)^\frac{2(q-1)}{2q-1}.
\end{equation}
Inserting (\ref{S.14}) and (\ref{S.15}) into (\ref{S.12}), Young's inequality
(since $\frac{1}{2q-1}<\frac{1}{q}$) yields
\begin{equation}\nonumber
\bigg(R^{-d}\int_{B_{2R}}|\eta u|^{2q}\bigg)^\frac{1}{q}\lesssim \bigg(R^{-d}\int_{B_{2R}}|u|\bigg)^{2}.
\end{equation}
Inserting this into (\ref{S.15}) and then into (\ref{S.12}) we get
\begin{equation}\nonumber
 R^{-d}\int_{B_{2R}}|R\nabla(\eta u)|^{2}\lesssim \bigg(R^{-d}\int_{B_{2R}}|u|\bigg)^{2},
\end{equation}
which by definition of $\eta$ turns into the first estimate in (\ref{S.17}).

\medskip

{\bf Step 3}.
We now address the somewhat subtle change necessary in Step 2 due to discreteness. The convenient continuum identity (\ref{nk.2}) can be substituted by the almost as convenient
discrete one
\begin{equation}\label{nk.3}
(\nabla(\eta^2u)a\nabla u)(b)=(\nabla(\eta u)a\nabla(\eta u))(b)-u(x)u(y)(a \nabla\eta)^2(b),
\end{equation}
where $x$ and $y$ denote the end points of the edge $b$.
We note that because of the diagonality of $a$, (\ref{nk.3}) reduces to the elementary identity
\begin{equation}\nonumber
(\eta^2u-\tilde\eta^2\tilde u)(u-\tilde u)-(\eta u-\tilde\eta\tilde u)^2=-u\tilde u(\eta-\tilde\eta)^2.
\end{equation}
Hence the discrete analogue of (\ref{S.12}) of Caccioppoli's estimate
is given by
\begin{equation}\nonumber
\bigg(\sum_{\E^d}(\nabla(\eta u))^2\bigg)^\frac{1}{2}\lesssim
\bigg(\sum_{b=[x,y]\in\E^d}u(x)u(y)(\nabla\eta(b))^2\bigg)^\frac{1}{2},
\end{equation}
which we rewrite in the dimensionless form of
\begin{equation}\label{nk.5}
\bigg(R^{-d}\sum_{\E^d}(R\nabla(\eta u))^2\bigg)^\frac{1}{2}\lesssim
\bigg(R^{-d}\sum_{b=[x,y]\in\E^d}u(x)u(y)(R\nabla\eta(b))^2\bigg)^\frac{1}{2}.
\end{equation}

\smallskip

We first turn to the l.\ h.\ s.\ of (\ref{nk.5}).
As in the continuum case, we appeal to the Poincar\'e-Sobolev estimate on $\mathbb{Z}^d$ applied
to the function $v=\eta u$ supported in $B_R$:
\begin{equation}\label{nk.4}
\bigg(R^{-d}\sum_{\Z^d}(\eta u)^{2q}\bigg)^\frac{1}{2q}\lesssim \bigg(R^{-d}\sum_{\E^d}(R\nabla(\eta
u))^2\bigg)^\frac{1}{2}.
\end{equation}
Again, such a discrete estimate can be derived from its continuum version (\ref{S.14})
by identifying $v=\eta u$ with a compactly supported finite element function on a triangulation subordinate to the lattice $\mathbb{Z}^d$. (Here this is easier than in Step~1 since $\eta u$ is supported in $B_{2R}$.)
On the r.\ h.\ s.\ of (\ref{nk.5}), we also proceed as in the continuum case and
apply H\"older's inequality:
\begin{multline}
R^{-d}\sum_{b=[x,y]\in\E^d}u(x)u(y)(R\nabla\eta(b))^2\\
\lesssim 
\bigg(R^{-d}\sum_{b=[x,y]\in\E^d}|u(x)|^q|u(y)|^q(R\nabla\eta(b))^{2\frac{2q-1}{q}}\bigg)^\frac{1}{2q-1}
\bigg(R^{-d}\sum_{B_{R+1}}                            |u|      \bigg)^\frac{2(q-1)}{2q-1}.\label{nk.a}
\end{multline}
Following the continuum case, we choose our cut-off function as
\begin{equation}\label{nk.b}
\eta(x)=\hat\zeta\Big(\frac{x}{R}\Big)^r\quad\mbox{with}\quad \frac{1}{r-1}+\frac{1}{q}=1,
\end{equation}
where we specify the mask $\hat\zeta$ to be $\hat\zeta(\hat x)=\max\{1-|\hat x|,0\}$. 
The discrete version of (\ref{neukamm}) reads
\begin{equation}\nonumber
|R\nabla\eta(b)|^\frac{2q-1}{q}\lesssim\max\{\eta(x),\eta(y)\}\quad\mbox{for}\;b=[x,y],
\end{equation}
which we use in form of
\begin{equation}\nonumber
|R\nabla\eta(b)|^\frac{2q-1}{q}\lesssim\min\{\eta(x),\eta(y)\}+R^{-r}.
\end{equation}
Hence we obtain (with help of the Cauchy-Schwarz inequality to separate $x$ and $y$)
\[
\sum_{b=[x,y]\in\E^d}|u(x)|^q|u(y)|^q(R\nabla\eta(b))^{2\frac{2q-1}{q}}
\lesssim\sum_{\Z^d}|\eta u|^{2q}+R^{-r}\sum_{B_{R+1}}|u|^{2q}.
\]
In view of \eqref{nk.b}, we can make $r>2qd$ by choosing $q>1$ very close to one (indeed $q\le \frac{2d+1}{2d}$ will do).
Hence by the discrete $\ell^{2q}$--$\ell^1$ inequality, the above turns into
\[
R^{-d}\sum_{b=[x,y]\in\E^d}|u(x)|^q|u(y)|^q(R\nabla\eta(b))^{2\frac{2q-1}{q}}
\lesssim R^{-d}\sum_{\Z^d}|\eta u|^{2q}+\bigg(R^{-d}\sum_{B_{R+1}}|u|\bigg)^{2q}.
\]
Inserting this inequality into (\ref{nk.a}) we obtain
\begin{multline*}
 R^{-d}\sum_{b=[x,y]\in\E^d}u(x)u(y)(R\nabla\eta(b))^2\\
\lesssim
\bigg(R^{-d}\sum_{\Z^d}|\eta u|^{2q}\bigg)^\frac{1}{2q-1}
\bigg(R^{-d}\sum_{B_{R+1}}                           |u|      \bigg)^\frac{2(q-1)}{2q-1}
+\bigg(R^{-d}\sum_{B_{R+1}}                        |u|      \bigg)^2.
\end{multline*}
We now deduce the desired result from this, (\ref{nk.4}) and (\ref{nk.5}).

\medskip

{\bf Step 4}. Proof of \eqref{quenched} in the continuum setting. By translation invariance, it is enough to prove (\ref{quenched}) for
$y=0$, that is,
\begin{equation}\label{S.23}
\int_{\frac{4}{3}R\le|x|\le\frac{8}{3}R}|\nabla_x G(x,0)|^2\;dx\lesssim R^{2-d}.
\end{equation}
Here comes the argument for (\ref{S.23}):
We first apply Step 1 to $u(x):=G(x,0)$, with $f(x) = \delta(x)$ and the ball $B_{4R}(0)$ around the origin.
Formally, we have that $\int_{\mathbb{R}^d}|f|=1$; we do not care about the lack of regularity of the Dirac distribution since this does not play a role in the discrete case. Hence estimate \eqref{S.20} translates to
\begin{equation}\label{S.21}
\int_{B_{4R}(0)}|u-\bar u|\lesssim R^2
\end{equation}
for some constant $\bar u\in\mathbb{R}$.
We then apply Step 2 to $u$ replaced by $u-\bar u$ and to the box $B_R(y)$ with an arbitrary point $y$ with $|y|=2R$
as center. Since $0\not\in B_R(y)$, the function $u$ satisfies (\ref{S.11}) in this ball, so that the result (\ref{S.17}) of
Step 2 turns into 
\begin{equation}\label{S.22}
\bigg(R^{-d}\int_{B_{R}(y)}|R\nabla u|^2\bigg)^\frac{1}{2}\lesssim R^{-d}\int_{B_{2R}(y)}|u-\bar u|
\stackrel{|y|=2R}{\le} R^{-d} \int_{B_{4R}(0)}|u-\bar u|.
\end{equation}
Inserting (\ref{S.22}) into (\ref{S.21}) yields
\begin{equation}\nonumber
\int_{B_R(y)}|\nabla u|^2\lesssim R^{2-d}\quad\mbox{for any}\;y\;\mbox{with}\;|y|=2R.
\end{equation}
Since the annulus $\{x:\frac{4}{3}R\le|x|\le\frac{8}{3}R\}$ can be covered by finitely many balls
of the form $\{B_R(y)\}_{|y|=2R}$ with a number only depending on $d$, we obtain (\ref{S.23}).

\smallskip

This step applies verbatim to the discrete setting. The only difficulty lies in the application of Step~1 to $G$, i.e.\ we need to
justify the integration by parts \eqref{int_by_parts}. Let $G_T(x,0)$ be the Green function with massive term $T>0$, i.e.\
$G_T$ is the solution to the weak (difference) equation
\[
 \frac{1}{T} \sum_{x\in\Z^d} \zeta(x) G_T(x,0) +  \sum_{e\in\E^d} \nabla \zeta(e) a(e) \nabla G_T(e,0) = \zeta(0),
\]
for all $\zeta$ with compact support. It holds that $G_T(x,0) = \int_0^\infty e^{-t/T}G(t,x,0)dt$ and therefore $\| G(t,\cdot,0) \|_{\ell^1} = 1$ yields $\| G_T \|_{\ell^1} = T$, i.e.\ $u_T:=G_T(\cdot,0)\in\ell^1(\Z^d)$. Since furthermore (in the notation of Step~1) $v_M \in \ell^\infty(\Z^d)$, the integration by
parts is valid and we may apply Step~1 to obtain \eqref{S.20} with $u$ replaced by $u_T$ uniformly in $T>0$. Since $\nabla G_T$ converges
point-wise to $\nabla G$ as $T\to\infty$, estimate \eqref{S.20} extends to $\nabla G$ by Fatou's lemma. In fact, the limit of $\nabla G_T$ may be taken as a definition of $\nabla G$ in the case of $d=2$.

\medskip

{\bf Step 5}. Proof of \eqref{quenched2} in its continuum version. By translation invariance, 
it is enough to prove it centered at the origin, that is,
\begin{equation}\label{S.25}
\int_{8R\le|y|\le 16R}\int_{|x|\le R}|\nabla\nabla G(x,y)|^2\;dxdy\lesssim 1,
\end{equation}
where we recall that $\nabla\nabla G$ denotes the second mixed derivative of $G$.
Here comes the argument for (\ref{S.25}). By symmetry of the Green function, we also have
\begin{equation}\nonumber
\nabla^*_y a(y)\nabla_y G(x,y)=0\quad\mbox{for}\;y\not=x,
\end{equation}
which we may differentiate w.\ r.\ t.\ $x$ to the effect of
\begin{equation}\label{S.26}
\nabla^*_y a(y)\nabla_y\nabla_x G(x,y)=0\quad\mbox{for}\;y\not=x.
\end{equation}
For fixed $x\in\R^d$ such that $|x|\le R$ and fixed $z$ with $|z|=12R$, we apply Step 2 to
the function $y\mapsto \nabla_x G(x,y)$ and the ball $B_{5R}(z)$. Since $|x-z|>10R$ 
and in view of (\ref{S.26}), this function is $a$-harmonic in $B_{10R}(z)$ so that we obtain from (\ref{S.17}) with $R$ replaced by $10R$
\begin{equation}\nonumber
\int_{|y-z|\le 5R}|\nabla\nabla G(x,y)|^2 \;dy \lesssim R^{-2} \int_{|y-z|\le10R}|\nabla_x G(x,y)|^2 \;dy.
\end{equation}
We integrate this estimate over $|x|\le R$:
\begin{equation}\nonumber
\int_{|y-z|\le 5R}\int_{|x|\le R}|\nabla\nabla G(x,y)|^2 \;dxdy \lesssim R^{-2} \int_{|y-z|\le 10R}\int_{|x|\le R}|\nabla_x G(x,y)|^2
\;dxdy.
\end{equation}
Since for $|x|\le R$ and $|z|=12R$ we have $R\le |x-y|\le 23R$ for all $|y-z|\le10R$, this turns into
%
\begin{multline*}
\int_{|y-z|\le 5R}\int_{|x|\le R}|\nabla\nabla G(x,y)|^2 \;dxdy\\ \lesssim
R^{-2} \int_{|y-z|\le 10R}\int_{R\le|x-y|\le 23R}|\nabla_x G(x,y)|^2\;dxdy.
\end{multline*}
Since the annulus $\{x:R\le|x-y|\le 23R\}$ can be covered by five dyadic annuli, we obtain from Step 4
\begin{equation}\nonumber
\int_{|y-z|\le 5R}\int_{|x|\le R}|\nabla\nabla G(x,y)|^2 \;dxdy \lesssim 1.
\end{equation}
Since the annulus $\{y:8R\le|y|\le 16R\}$ can be covered by finitely many balls
of the form $\{B_{5R}(z)\}_{|z|=12R}$ with a number only depending on $d$, we obtain (\ref{S.25}). All estimates in this step hold verbatim
in the discrete setting.

\section{Proof of Proposition \ref{prop:DD_ell}}\label{sec:DD_ell}

The proof we give differs from the original one of Delmotte and Deuschel in that
it relies on quenched regularity for the elliptic Green function, cf.\ Lemma \ref{lem:quenched}, 
rather than the parabolic one. Again, $\lesssim$ means
up to a generic constant that only depends on $d$
and $\lambda$. We shall prove all statements in their continuum version first and then indicate the changes for the discrete setting. Again, we shall not worry about the finiteness and regularity of the continuum expressions.

\medskip

{\bf Step 1}. We claim that by stationarity of $\langle\cdot\rangle$, we have for any
shift vector $z\in\mathbb{R}^d$ and any exponent $p$
\begin{eqnarray}
\langle|\nabla_xG(x,y)|^p\rangle&=&\langle|\nabla_xG(x+z,y+z)|^p\rangle,\label{S.35}\\
\langle|\nabla\nabla G(x,y)|^p\rangle&=&\langle|\nabla\nabla G(x+z,y+z)|^p\rangle,\label{S.36}
\end{eqnarray}
In particular, for the continuum version of Proposition \ref{prop:DD_ell}, it is enough to show
\begin{align}
\langle|\nabla_x G(0,y)|^2\rangle^\frac{1}{2}&\lesssim|y|^{1-d},\label{S.45}\\
\langle|\nabla\nabla G(0,y)|\rangle&\lesssim|y|^{-d}.\label{S.46}
\end{align}
Indeed, by uniqueness of the Green function
we have for any shift vector $z\in\mathbb{R}^d$
\begin{equation}\nonumber
G(a(\cdot+z);x,y)=G(a;x+z,y+z),
\end{equation}
which we may differentiate and take to the $p$-th power to obtain
\begin{align*}
|\nabla_xG(a(\cdot+z);x,y)|^p&=|\nabla_xG(a;x+z,y+z)|^p,\\
|\nabla\nabla G(a(\cdot+z);x,y)|^p&=|\nabla\nabla G(a;x+z,y+z)|^p.
\end{align*}
Hence by stationarity of $\langle\cdot\rangle$, this implies (\ref{S.35}) and (\ref{S.36}). The same reduction can be made in the discrete
setting.

\medskip

{\bf Step 2}. We claim that by Lemma \ref{lem:quenched} we have for any radius $R$
\begin{equation}\label{S.37}
\bigg\langle R^{-d}\int_{y:R\le|y|\le2R} \Big( |\nabla\nabla G(0,y)|^2+R^{-2}|\nabla_xG(0,y)|^2 \Big) \;dy \bigg\rangle^\frac{1}{2}
\lesssim R^{-d}.
\end{equation}
%
%
%
For the first summand, we appeal to \eqref{quenched2}, which we use in its continuum version, i.e.~\eqref{S.25},
\begin{equation}\nonumber
\bigg(R^{-2d}\int_{8R\le|y|\le 16R}\int_{|x|\le R}|\nabla\nabla G(x,y)|^2\;dxdy\bigg)^\frac{1}{2}\lesssim R^{-d}.
\end{equation}
Adding to it its version with $R$ replaced by $2R$ we obtain a similar statement with a thicker annulus:
\begin{equation}\nonumber
\bigg(R^{-2d}\int_{|x|\le R}\int_{8R\le|y|\le 32R}|\nabla\nabla G(x,y)|^2\;dydx\bigg)^\frac{1}{2}\lesssim R^{-d}.
\end{equation}
Taking the square expectation yields
\begin{equation}\nonumber
\bigg(R^{-2d}\int_{|x|\le R}\int_{8R\le|y|\le 32R}\langle|\nabla\nabla G(x,y)|^2\rangle \;dydx\bigg)^\frac{1}{2}
\lesssim R^{-d},
\end{equation}
which by (\ref{S.36}) takes the form of
\begin{equation}\nonumber
\bigg(R^{-2d}\int_{|x|\le R}\int_{8R\le|y|\le 32R}\langle|\nabla\nabla G(0,y-x)|^2\rangle \;dydx\bigg)^\frac{1}{2}
\lesssim R^{-d}.
\end{equation}
Since for any $x$ with $|x|\le R$, $y'=y-x$ covers the annulus $\{y':12R\le|y'|\le24R\}$ if
$y$ runs through the annulus $\{y:8R\le|y|\le 32R\}$, this implies as desired
\begin{equation}\nonumber
\bigg(R^{-d}\int_{12R\le|y'|\le 24R}\langle|\nabla\nabla G(0,y')|^2\rangle \;dy'\bigg)^\frac{1}{2}
\lesssim R^{-d},
\end{equation}
i.e.~\eqref{S.37} with $R$ replaced by $12R$.
For the second estimate, we just use stationarity in form of $\langle |\nabla_x G(0,y)|^2 \rangle \stackrel{\eqref{S.35}}{=} \langle |\nabla_y G(-y,0)|^2 \rangle$, to obtain
\[
 \bigg\langle \int_{y:R\le|y|\le2R} |\nabla_xG(0,y)|^2 \;dy \bigg\rangle = \bigg\langle \int_{y:R\le|y|\le2R} |\nabla_y G(y,0)|^2 \;dy \bigg\rangle.
\]
Hence by \eqref{quenched} in its continuum version, i.e.~\eqref{S.23}, the second summand in (\ref{S.37}) is under control. All estimates remain valid in the discrete case.

\medskip

{\bf Step 3}. Consider the $a$-dependent functions $u=u(a;x)$, $f=f(a;x)$ and the vector
field $g=g(a;x)$ related by
\begin{equation}\label{S.30}
\nabla^*a\nabla u=\nabla^*g+f\quad\mbox{in}\;\R^d.
\end{equation}
Suppose that $f$ and $g$ are supported on an annulus of radius $R$:
\begin{equation}\label{S.31}
f(x)=0,\; g(x)=0\quad\text{unless}\;R\le|x|\le2R.
\end{equation}
Then we claim
\begin{align}
\langle|\nabla u(0)|^2\rangle^\frac{1}{2}&\lesssim\sup_{a\in\Omega}\bigg(R^{-d}\int_{\R^d} \big( |g|^2 + R^2 |f|^2 \big)
\bigg)^\frac{1}{2},\label{S.32}\\
\langle|\nabla u(0)|  \rangle            &\lesssim\Big\langle R^{-d}\int_{\R^d} \big( |g|^2 + R^2 |f|^2 \big)
\Big\rangle^\frac{1}{2}.\label{S.33}
\end{align}
To prove \eqref{S.32} and \eqref{S.33}, we start by noting that (\ref{S.30}) yields the representation formula
\begin{equation}\nonumber
u(x)=\int_{\R^d} G(x,y)(\nabla^*g+f)(y) \;dy,
\end{equation}
which we use in form of
\begin{equation}\nonumber
\nabla u(0)= \int_{\R^d} \nabla\nabla G(0,y) g(y) \;dy + \int_{\R^d} \nabla_x G(0,y) f(y) \;dy.
\end{equation}
By the  Cauchy-Schwarz inequality in space and the support assumption, this yields
\begin{multline*}
|\nabla u(0)|\le\bigg(\int_{R\le|y|\le2R}|\nabla\nabla G(0,y)|^2 \;dy \int_{R\le|y|\le2R} |g(y)|^2 \;dy \bigg)^{\frac{1}{2}}\\
+ \bigg( \int_{R\le|y|\le2R} |\nabla_x G(0,y)|^2 \;dy \int_{R\le|y|\le2R} |f(y)|^2 \;dy \bigg)^\frac{1}{2}.
\end{multline*}
This implies by the Cauchy-Schwarz in probability:
\begin{align}
\langle|\nabla u(0)|^2\rangle^\frac{1}{2}&\le\Lambda\sup_{a\in\Omega}\bigg(R^{-d}\int_{\R^d} \big( |g|^2+R^2|f|^2 \big)
\bigg)^\frac{1}{2},\nonumber\\
\langle|\nabla u(0)|  \rangle            &\le\Lambda\Big\langle R^{-d}\int_{\R^d} \big( |g|^2+R^2 |f|^2 \big)
\Big\rangle^\frac{1}{2},\nonumber
\end{align}
where we have set for abbreviation
\begin{equation}\nonumber
\Lambda:=\bigg\langle R^d \int_{R\le|y|\le2R} \big( |\nabla\nabla G(0,y)|^2 + R^{-2} |\nabla_x G(0,y)|^2 \big) \;dy
\bigg\rangle^\frac{1}{2}.
\end{equation}
Hence in order to obtain (\ref{S.32}) and (\ref{S.33}), we need $\Lambda\lesssim 1$
which is just the statement of Step 2. All estimates in this step carry over to the discrete setting.

\medskip

{\bf Step 4}. Consider an $a$-dependent functions $u=u(a;x)$ satisfying
\begin{equation}\label{S.42}
\nabla^*a\nabla u=0\quad\mbox{in}\;B_{2R}
\end{equation}
on the ball $B_{2R}$ of radius $2R$ around the origin. Then we claim
\begin{align}
\langle|\nabla u(0)|^2\rangle^\frac{1}{2}&\lesssim\sup_{a\in\Omega}\bigg(R^{-d}\int_{B_{2R}}|\nabla u|^2\bigg)^\frac{1}{2},
\label{S.40}\\
\langle|\nabla u(0)|  \rangle            &\lesssim\Big\langle R^{-d}\int_{B_{2R}}|\nabla u|^2\Big\rangle^\frac{1}{2}.
\label{S.41}
\end{align}
To see this, pick a cut-off function $\eta$ for $B_R$ in $B_{2R}$ such that $|\nabla\eta| \lesssim R^{-1}$ and set
\begin{equation}\nonumber
v:=\eta(u-\bar u),\quad\mbox{where}\;\bar u\;\mbox{denotes average of $u$ on $B_{2R}$}.
\end{equation}
Equation \eqref{S.42} yields
\begin{equation}\label{u_harmonic_f_g}
 \nabla^* a \nabla v = \nabla^*g + f\quad\text{with}\quad g:=(u-\bar u)a\nabla\eta\quad\text{and}\quad f:=-\nabla\eta\cdot a\nabla u.
\end{equation}
By choice of $\eta$, the functions $g$ and $f$ satisfy the support condition (\ref{S.31}) and we have
\begin{equation}\nonumber
\int_{\R^d} \big( |g|^2 + R^2 |f|^2 \big) \lesssim \int_{B_{2R}} \big( R^{-2} |u-\bar u|^2 + |\nabla u|^2 \big).
\end{equation}
Hence the Poincar\'{e} inequality on $B_{2R}$ applied to the first term on the r.\ h.\ s.\ yields
\[
 \int_{\R^d} \big( |g|^2 + R^2 |f|^2 \big) \lesssim \int_{B_{2R}}|\nabla u|^2
\]
Thus (\ref{S.40}) and (\ref{S.41}) follow from (\ref{S.32}) and (\ref{S.33}).

\smallskip

In the discrete setting, we just need to check that \eqref{u_harmonic_f_g} remains valid even in the absence of the (continuum) Leibniz
rule.
Indeed, for any functions $a$, $v$ and $\eta$ we have that
\begin{multline*}
 \nabla^* \big( a \nabla( \eta v)\big)(x) + \sum_{j=1}^d \nabla\eta([x,x+e_j]) a([x,x+e_j]) \nabla v([x,x+e_j])\\
 = \nabla^*(a v \nabla \eta)(x) + \eta(x) (\nabla^* a \nabla v)(x)
\end{multline*}
for all $x\in\Z^d$, where we have set $(a v \nabla \eta)([x,x+e_j]) = a([x,x+e_j]) v(x) \nabla \eta([x,x+e_j])$ for any edge $[x,x+e_j]\in\E^d$.
 By definition of $\nabla$ and $\nabla^*$, this follows from the elementary identity
\begin{align*}
 &a([x-e_j,x]) \big( v(x) \eta(x) - v(x-e_j) \eta(x-e_j) \big)\\
 &\qquad- a([x,x+e_j]) \big( v(x+e_j) \eta(x+e_j) - v(x) \eta(x) \big)\\
 &\qquad+ \big( \eta(x+e_j) - \eta(x) \big) a([x,x+e_j]) \big( v(x+e_j) - v(x) \big)\\
 &= a([x-e_j,x]) v(x-e_j) \big( \eta(x) - \eta(x-e_j) \big) - a([x,x+e_j]) v(x) \big(\eta(x+e_j) - \eta(x) \big)\\
 &\qquad+ \eta(x) a([x-e_j,x]) \big( v(x) - v(x-e_j) \big) - \eta(x) a([x,x+e_j]) \big( v(x+e_j) - v(x) \big).
\end{align*}
If $\eta$ is a cut-off function for $B_{R+1}$ in $B_{2R-1}$, this shows that \eqref{u_harmonic_f_g} is valid also in the discrete case with the desired support condition on $f$ and $g$.

\medskip

{\bf Step 5}. Conclusion, that is, proof of (\ref{S.45}) and (\ref{S.46}).
We fix $y\in\mathbb{R}^d\setminus \{0\}$ and apply Step 4 to
$u(x)=G(x,y)$ and $R=\frac{1}{6}|y|$.
From (\ref{S.40}) we obtain
\begin{equation}\nonumber
\langle|\nabla_x G(0,y)|^2\rangle^\frac{1}{2}\lesssim\sup_{a\in\Omega}\bigg(|y|^{-d}\int_{B_{\frac{1}{3}|y|}}
|\nabla_x G(x,y)|^2 \;dx \bigg)^\frac{1}{2}.
\end{equation}
Since
\begin{equation}\label{ball_ann}
 B_{\frac{1}{3}|y|}\subset\Big\{x\in\R^d:\frac{2}{3}|y|\le|x-y|\le\frac{4}{3}|y|\Big\},
\end{equation}
we obtain by (\ref{quenched}) in its continuum version \eqref{S.23} with $R=\frac{1}{2}|y|$ as desired (\ref{S.45}), i.e.\ we have that
\begin{equation}\nonumber
\langle|\nabla_x G(0,y)|^2\rangle^\frac{1}{2}\lesssim |y|^{1-d}.
\end{equation}
We now apply Step 4 to the $a$-harmonic function $u(x)=\nabla_y G(x,y)$ with $R=\frac{1}{6}|y|$ and obtain from (\ref{S.41}) that
\begin{equation}\nonumber
\langle|\nabla\nabla G(0,y)|\rangle\lesssim
\bigg\langle|y|^{-d}\int_{B_{\frac{1}{3}|y|}}|\nabla\nabla G(x,y)|^2 \;dx\bigg\rangle^\frac{1}{2}.
\end{equation}
The inclusion \eqref{ball_ann} yields
\begin{equation}\nonumber
\langle|\nabla\nabla G(0,y)|\rangle\lesssim
\bigg\langle|y|^{-d}\int_{\frac{2}{3}|y|\le|x-y|\le\frac{4}{3}|y|}|\nabla\nabla
G(x,y)|^2 \;dx \bigg\rangle^\frac{1}{2}.
\end{equation}
To conclude, we want to apply \eqref{S.37}. In view of (\ref{S.36}) and the symmetry of $\nabla\nabla G$, we rewrite \eqref{S.37} as
\begin{equation}\nonumber
\bigg\langle R^{-d}\int_{R\le|x-y|\le2R}|\nabla\nabla G(x,y)|^2 \;dx \bigg\rangle^\frac{1}{2}
\lesssim R^{-d}.
\end{equation}
Letting $R=\frac{2}{3}|y|$ yields (\ref{S.46}). This step carries over verbatim to the discrete setting if $|y|$ is large enough (which we used a few times in previous steps). The
conclusion for the finitely many small $y$ follows from the quenched bounds on $G$, i.e.\ \eqref{nash}.


{\small
\begin{thebibliography}{999}
%

\bibitem{Aro} {\it D.\ G.\ Aronson}: {Bounds for the fundamental solution of a parabolic equation}. Bull.\ Amer.\ Math.\ Soc.\ {\bf 73} (1967), 890--896. Zbl 0153.42002, MR0217444.
\bibitem{Del} {\it T.\ Delmotte}: Parabolic {H}arnack inequality and estimates of {M}arkov chains on graphs. Rev.\ Mat.\ Iberoamericana {\bf 15} (1999), no.\ 1, 181--232. Zbl 0922.60060, MR1681641.
\bibitem{DD} {\it T.\ Delmotte and J.-D.\ Deuschel}: {On estimating the derivatives of
symmetric diffusions in stationary random environments, with applications to the $\nabla\phi$ interface model}. Probab.\ Theory Relat.\ Fields {\bf 133} (2005), 358--390. Zbl 1083.60082, MR2198017.
\bibitem{LNO} {\it A.\ Lamacz, S.\ Neukamm and F.\ Otto}: {Moment bounds for the corrector in stochastic homogenization of a percolation model}. WIAS Preprint No.\ 1863 (2013).
\bibitem{MO} {\it D.\ Marahrens and F.\ Otto}: {Annealed estimates on the Green function}. \verb#arXiv:1304.4408# (2013).
\bibitem{Nash} {\it J.\ Nash}: Continuity of solutions of parabolic and elliptic equations. American J.\ Math.\ {\bf 80} (1958), 931--954. Zbl 0096.06902, MR0100158.
\end{thebibliography}}
\end{document}